\documentclass[
a4paper,				
english
]{extarticle}

\newif\ifshow
\showtrue   

\newif\ifmarkchanges
\markchangestrue   

\makeatletter
\title{Computing with Orientation-Dependent Vector Fields in Crystallographic Texture Analysis}\let\Title\@title
\date{\today}\let\Date\@date
\makeatother
\newcommand{\keywords}{rotation group SO(3),crystallographic texture analysis,harmonic analysis,Wigner-D functions,generalized spherical harmonics,rotation vectors,vector fields on SO(3),gradient,divergence,curl,crystal symmetry,MTEX}

\usepackage{header}
\usepackage{headerFormulas}

\begin{document}

\author{Ralf
  Hielscher\orcidlink{0000-0002-6342-1799}%
  \thanks{ralf.hielscher@math.tu-freiberg.de}~,~%
  Erik Wünsche \orcidlink{0000-0001-5178-6654}%
  \thanks{erik.wuensche@math.tu-freiberg.de}}

\maketitle
\begin{abstract}
  We present a computational framework for representing orientation-dependent vector fields.
  Building on~\cite{Hielscher2026b}, we translate the differential-geometric formulation into crystallographic notation and focus on its computational use in crystallographic texture analysis.

  In crystallographic applications, small orientation changes are commonly described by rotation vectors, spin tensors, or infinitesimal misorientations.
  By contrast, vector fields in harmonic computations are usually represented by components with respect to Euler-angle derivatives, following the classical formulation of generalized spherical harmonics.
  This representation is convenient but numerically unfavorable.
  We instead use rotation-vector components, extending a familiar description of individual orientation changes to orientation-dependent vector fields.

  In this representation, vector fields are described by three scalar component functions, which can be expanded into generalized spherical harmonics.
  We derive frequency-domain formulas for gradient, divergence, and curl, enabling their evaluation directly on harmonic coefficients.
  In particular, the gradient of a band-limited scalar function remains band-limited with the same bandwidth.
\end{abstract}

\tableofcontents

\section{Introduction}\label{sec:Introduction}

The study of crystallographic microstructures and their deformation is largely based on the analysis of orientation-valued data.
More precisely, in a polycrystalline material, each crystallite is described by its crystallographic orientation relative to a fixed specimen coordinate system.
This orientation describes the change from the crystal frame, fixed to the lattice, to the specimen frame.
Such a change of frame is described by a rotation, and since both frames may carry symmetries, an orientation is a rotation together with these symmetries.
The set of all rotations is called the \emph{rotation group} and denoted by $\SO3$. For simplicity, we also refer to $\SO3$ as the \emph{orientation space} in the following.
Crystallographic texture is then described as a distribution of orientations on $\SO3$.
The mathematical foundations of orientations and their symmetries, their different representations, and orientation distributions have been developed in detail by Bunge~\cite{Bunge1982} and Morawiec~\cite{Morawiec2004}.

In many computations in crystallographic texture analysis, derivatives or rates of change of orientation-valued data occur naturally. Important examples include the following:
\begin{itemize}
  \item In time-dependent texture evolution~\cite{Morawiec2004,Clement1979,Bunge1984}, the change of an orientation distribution function (ODF) is governed by the reorientation of crystallites, which is described by an orientation-dependent spin on $\SO3$.
  \item Spatially resolved orientation data lead to spatial gradients of orientation maps. Such gradients are used, for instance, in the analysis of lattice curvature~\cite{Sun2000}, geometrically necessary dislocations~\cite{Sun2000,Pantleon2008}, and Nye tensor computations~\cite{Nye1953,Sun2000,Pantleon2008}.
  Although these derivatives are taken with respect to physical space, they describe local changes of crystallographic orientation.
  \item Many optimization tasks involving orientation-dependent quantities require gradients with respect to orientation. Examples include the optimization of orientation-dependent material properties~\cite{Gu2018}, the training of orientation-dependent neural network models~\cite{Ding2020}, and the construction of representative orientation sets by gradient-based methods~\cite{Graef2013,Knezevic2024}.
\end{itemize}

In all these examples, vector fields on $\SO3$ arise naturally. Such a vector field is a function that assigns to each orientation $\mat R\in\SO3$ an object that describes an orientation-dependent infinitesimal change of orientation.
There are several equivalent ways to describe changes of orientation, for example by spin tensors, misorientations, or left- and right-trivialized rotation vectors.
In the present paper, we use a common geometric terminology and regard infinitesimal changes of orientation as tangent vectors on $\SO3$.

The efficient representation and evaluation of orientation-dependent vector fields is an essential component of computational texture analysis.
As a starting point, we first recall the simpler case of scalar orientation-dependent functions.
Important examples are orientation distribution functions (ODFs) and orientation-dependent material properties, such as elastic moduli~\cite{Mainprice2011}, yield stresses~\cite{Taylor1938,Bishop1951a}, or corrosion-related quantities~\cite{Shin2012}.
Such functions $f\colon\SO3\to\IR$ are most often represented by expansions into generalized spherical harmonics, i.e.
\begin{equation*}
  f(\mat R) = \sum_{n=0}^{\infty}\sum_{k,k'=-n}^n \fhat{n}{k}{k'} \, \WignerD{n}{k}{k'}(\mat R),
\end{equation*}
where $\fhat{n}{k}{k'}  = \scp[\big]{f}{\WignerD{n}{k}{k'}}_{\L2SO3}$ are the harmonic coefficients of $f$.
This representation was established as a standard tool in crystallographic texture analysis by Bunge~\cite{Bunge1982}.
Instead of storing function values or evaluating the function from a pointwise model, the function is represented by its harmonic coefficients.
This gives a compact spectral representation and enables efficient evaluation and computation of basic quantities, such as the texture index~\cite{Bunge1982}, volume portions~\cite{Zuo1993,Cho2004}, cross correlations~\cite{Hielscher2019,Lenthe2020}, and pole figures~\cite{Hielscher2008}, within the harmonic framework.
Moreover, fast \(\SO3\)-Fourier algorithms are available for evaluating such expansions and for computing their coefficients~\cite{Risbo1996,Kostelec2008,Potts2009,Hielscher2026}.
In practice, the series is truncated at a finite bandwidth, and the resulting truncation error can generally be reduced by increasing this bandwidth.

The goal of this paper is to discuss a computationally efficient way to extend the harmonic approximation scheme from scalar orientation-dependent functions to vector fields on the orientation space.
The starting point, explained in \Cref{sec:TangentVectors}, is to describe every infinitesimal change of orientation consistently by either its left or its right rotation vector.
Assigning such a rotation vector to every orientation yields a vector field $\vec f\colon\SO3\to\IR^{3}$ whose three component functions are ordinary scalar orientation-dependent functions, as discussed in \Cref{sec:VectorFields}.
Crystal and specimen symmetries act on these components differently than on scalar functions, which is the subject of \Cref{sec:symmetries}.
Since the components are scalar functions, they can be expanded into generalized spherical harmonics and treated by the same harmonic tools as an ODF.
This leads to the harmonic vector fields of \Cref{sec:HarmonicVectorFields}.
Finally, \Cref{sec:GradientDivCurl} shows that gradient, divergence, and curl can then be computed exactly and efficiently, directly from the harmonic coefficients.

Our aim is to establish this as the standard representation for vector fields on $\SO3$, thereby complementing the computational framework for crystallographic textures described by Morawiec~\cite{Morawiec2004} towards tangent vector fields and differential operators.
To this end, we build on~\cite{Hielscher2026b}, where vector fields and differential operators on $\SO3$ are developed rigorously in the language of differential geometry.
The main contribution of this tutorial is not to derive new results, but to translate this framework into the language of crystallographic texture analysis, i.e.\ into orientations and rotation vectors.
We focus on the ideas and formulas needed for computations and refer to~\cite{Hielscher2026b} for the proofs.

In the literature, orientation-dependent vector fields are often represented by their components with respect to derivatives along the Euler angles~\cite{Bunge1986,Kalidindi2005,Knezevic2024}.
This choice is natural, since generalized spherical harmonics are defined in terms of Euler angles, but it is numerically unfavorable.
We therefore recommend representing vector fields by their rotation-vector components rather than by their Euler-angle derivatives.

The framework investigated in this paper has been implemented in \texttt{MTEX}~\cite{MTEX} since version~6.0, and is thus the standard way in which \texttt{MTEX} represents and computes with orientation-dependent vector fields.


\section{Orientation Gradients and Their Representation}\label{sec:TangentVectors}

As mentioned above, several descriptions of small changes of orientation are used in crystallographic texture analysis, including spin tensors, misorientations, and left and right rotation vectors.
Let $\mat R_{0}\in\SO3$ be an initial orientation and let $\mat R_{1}\in\SO3$ be a nearby target orientation.
The corresponding relative rotations can be written as $\Delta_{L} = \mat R_{1}\mat R_{0}^{\top}$ and $\Delta_{R} = \mat R_{0}^{\top}\mat R_{1}$, so that
\[ \mat R_{1} = \Delta_{L} \, \mat R_{0} \qquad \text{and} \qquad \mat R_{1} = \mat R_{0} \, \Delta_{R}. \]
Here, $\Delta_{L}$ describes the change of orientation by multiplication from the left, whereas $\Delta_{R}$ describes the same change by multiplication from the right.

Since $\Delta_{L}$ and $\Delta_{R}$ are again rotations, they can be represented by their rotation axes and rotation angles.
This gives the left and right rotation vectors
\[ \vec\rho_{L} = \omega(\Delta_{L}) \operatorname{axis}(\Delta_{L}), \qquad \vec\rho_{R} = \omega(\Delta_{R}) \operatorname{axis}(\Delta_{R}), \]
where the rotation axes are normalized to unit length.

The two vectors describe the same change of orientation, but with respect to different coordinate systems.
The vector $\vec\rho_{R}$ describes the change in crystal coordinates, whereas $\vec\rho_{L}$ describes it in specimen coordinates.
They are related by
\begin{equation}\label{eq:ChangeLeftRightTangentVectors}
  \vec\rho_{L} = \mat R_{0} \, \vec\rho_{R}.
\end{equation}

Mathematically, a small change of orientation at $\mat R_{0}$ is called a \emph{tangent vector} to $\SO3$ at $\mat R_{0}$, see~\cite{Hielscher2026b}.
Misorientations, rotation vectors, and spin tensors are different ways to describe this one object in computations.
In differential geometry, tangent vectors to the orientation space are represented by skew-symmetric matrices
\[ \mathfrak{so}(3) \coloneqq \set{\mat S\in\IR^{3\times3}}{\mat S = -\mat S^{\top}}. \]
More precisely, a rotation vector
\[ \vec\rho = \begin{pmatrix}  x \\ y \\ z \end{pmatrix} \]
is linked to the skew-symmetric matrix
\[ \widehat{\rho} = \begin{pmatrix} 0 & -z & y \\ z & 0 & -x \\ -y & x & 0 \end{pmatrix}. \]
For sufficiently small orientation changes, the relative rotations satisfy $\Delta_{L} \approx \id + \widehat{\rho}_{L}$ and $\Delta_{R} \approx \id + \widehat{\rho}_{R}$.
Consequently,
\begin{equation*}
  \mat R_{1} - \mat R_{0} = \Delta_{L} \, \mat R_{0} - \mat R_{0} \approx \widehat{\rho}_{L} \, \mat R_{0}
  \qquad \text{and} \qquad
  \mat R_{1} - \mat R_{0} = \mat R_{0} \, \Delta_{R} - \mat R_{0} \approx \mat R_{0} \, \widehat{\rho}_{R} .
\end{equation*}
Thus, the same tangent vector can be described either by a rotation vector, by the corresponding skew-symmetric matrix, or by the induced first-order change of the orientation matrix.
In compact form, this correspondence can be summarized as
\begin{equation*}
  \begin{alignedat}{3}
    \vec\rho_{L} = \begin{pmatrix} x_{L}\\y_{L}\\z_{L} \end{pmatrix} \quad &\Longleftrightarrow \quad \widehat\rho_{L} = \begin{pmatrix} 0 & -z_{L} & y_{L} \\ z_{L} & 0 & -x_{L}\\ -y_{L} & x_{L} & 0 \end{pmatrix} \quad &\Longleftrightarrow \quad \mat R_{1}-\mat R_{0} \approx \widehat\rho_{L} \, \mat R_{0}\\
    \vec\rho_{R} = \begin{pmatrix} x_{R}\\y_{R}\\z_{R} \end{pmatrix} \quad &\Longleftrightarrow \quad \widehat\rho_{R} = \begin{pmatrix} 0 & -z_{R} & y_{R} \\ z_{R} & 0 & -x_{R}\\ -y_{R} & x_{R} & 0 \end{pmatrix} \quad &\Longleftrightarrow \quad \mat R_{1}-\mat R_{0} \approx \mat R_{0} \, \widehat\rho_{R}
  \end{alignedat}
\end{equation*}

It is important to note that a rotation vector, or equivalently the corresponding skew-symmetric matrix, represents an orientation gradient only at a specific orientation.
This is visible on the right-hand side of the two equivalence relations above: the actual change of the orientation matrix is obtained by applying the skew-symmetric matrix either from the left or from the right to the current orientation.
Thus, an infinitesimal reorientation is not described by a rotation vector alone, but by the pair consisting of the current orientation and the corresponding rotation vector.
In the left representation, this pair is $(\mat R_{0},\vec\rho_{L})$, whereas in the right representation it is $(\mat R_{0},\vec\rho_{R})$.

In the following, we want to represent vector fields with respect to a basis and define the gradient as the direction of steepest increase.
The former requires a vector space, in which every tangent vector can be written as a linear combination of three basis vectors.
The latter requires a notion of length and angle, i.e., an inner product.

For a fixed orientation $\mat R\in\SO3$, the set of all tangent vectors attached to $\mat R$ forms a vector space, which is called the tangent space at $\mat R$.
Using the rotation-vector representation introduced above, we equip this space with the inner product
\[ \scp{\vec\rho_{L}}{\vec\eta_{L}}_{2} = \frac12 \scp{\widehat\rho_{L}}{\widehat\eta_{L}}_{\mathrm F} = \frac12 \scp{\widehat\rho_{L}\, \mat R}{\widehat\eta_{L} \, \mat R}_{\mathrm F}.  \]
Here, \(\scp{\cdot}{\cdot}_{2}\) denotes the Euclidean inner product of rotation vectors, and \(\scp{\cdot}{\cdot}_{\mathrm F}\) denotes the Frobenius inner product of matrices.
Equivalently, in the right representation, one has
\[ \scp{\vec\rho_{R}}{\vec\eta_{R}}_{2} = \frac12 \scp{\widehat\rho_{R}}{\widehat\eta_{R}}_{\mathrm F} = \frac12 \scp{\mat R \, \widehat\rho_{R}}{\mat R \, \widehat\eta_{R}}_{\mathrm F}.  \]


\section{Orientation-Dependent Vector Fields}\label{sec:VectorFields}

We now consider mappings that assign to every orientation $\mat R\in\SO3$ a tangent vector at $\mat R$.
Such mappings are called tangent vector fields, or simply \emph{vector fields}.
A basic example is the gradient of an ODF.
It assigns to each orientation $\mat R$ the tangent vector that points in the direction of steepest increase of the ODF at \(\mat R\).

Since we have already introduced several descriptions of tangent vectors, there are several descriptions of a vector field. It may assign to each orientation \(\mat R\in\SO3\)
\begin{enumerate}[label=(\alph*)]
  \item a geometric tangent vector, i.e., an element of \(\mathfrak{so}(3)\cdot\mat R = \mat R\cdot\mathfrak{so}(3)\),
  \item a skew-symmetric matrix,
  \item a rotation vector in \(\IR^3\).
\end{enumerate}
From a geometric point of view, description~(a) is the most intrinsic one. However, since the three descriptions are equivalent, from a computational point of view we work with description~(c), i.e., we regard vector fields as functions $\SO3\to\IR^{3}$.

To arrive at this description, we first consider how general vector fields can be represented on the whole orientation space $\SO3$.
To this end, we return to the concept of a basis of a vector space.
We choose three smooth vector fields $\vec X_{1}$, $\vec X_{2}$, $\vec X_{3}$ such that for every orientation $\mat R\in\SO3$, the three corresponding tangent vectors
\[ \vec X_{1}(\mat R),\quad \vec X_{2}(\mat R),\quad \vec X_{3}(\mat R) \]
form a basis of all tangent vectors attached to $\mat R$, i.e., of the tangent space at $\mat R$.
In the differential-geometric terminology used in~\cite{Hielscher2026b}, such a smoothly varying choice of bases is called a frame.
Here, we continue to call it a basis.

With respect to this basis, every vector field $\vec f$ is represented by its component functions $f_{1},f_{2},f_{3}\colon\SO3\to\IR$ such that
\begin{equation}\label{eq:LinKombiVectorField}
  \vec f(\mat R) = f_{1}(\mat R) \, \vec X_{1}(\mat R) + f_{2}(\mat R) \, \vec X_{2}(\mat R) + f_{3}(\mat R) \, \vec X_{3}(\mat R),
\end{equation}
where $f_{1},f_{2},f_{3}\colon\SO3\to\IR$ are scalar orientation-dependent functions.

\paragraph{Which properties should the basis have?}
The choice of this basis is crucial for numerical computations.
The three vector fields should be geometrically meaningful, easy to evaluate, and well behaved on the whole orientation space.
In particular, we require the following properties:

\begin{enumerate}[label=\arabic*~]
  \item The vector fields $\vec X_{1},\vec X_{2},\vec X_{3}$ are smooth.
  \item Each vector field has unit length, i.e., $\abs{\vec X_i(\mat R)} = 1$ for all $\mat R\in\SO3$ and $i=1,2,3$.
  \item The basis degenerates nowhere.
        That is, for every orientation $\mat R\in\SO3$, every tangent vector attached to $\mat R$ can be represented uniquely as a linear combination
        \[ a \, \vec X_{1}(\mat R) + b \, \vec X_{2}(\mat R) + c \, \vec X_{3}(\mat R) \]
        with $a,b,c\in\IR$.
  \item[(4)] For every fixed orientation $\mat R\in\SO3$, the three tangent vectors $\vec X_{1}(\mat R),\vec X_{2}(\mat R),\vec X_{3}(\mat R)$ are mutually orthogonal.
\end{enumerate}
Property~(4) is not strictly necessary.
Together with the normalization in property~2, however, it yields an orthonormal basis at every orientation.
This simplifies the representation of vector fields and improves the numerical stability of the resulting computations.

\paragraph{How is the basis usually chosen?}
A computationally natural choice, frequently used in crystallographic texture analysis~\cite{Bunge1986,Kalidindi2005,Knezevic2024}, is obtained from infinitesimal orientation changes along the Euler angles $\phi_{1},\Phi,\phi_{2}$.
The reason is that the generalized spherical harmonics, in which ODFs are classically expanded~\cite{Bunge1982}, are defined in terms of Euler angles.

Let $\mat R(\phi_1,\Phi,\phi_2)$ denote the Euler-angle representation of an orientation.
Then the corresponding three vector fields are given by
\begin{equation*}
  \vec X_{1}(\mat R) = \left.\frac{\d{}}{\d{t}}\mat R(\phi_{1}+t,\Phi,\phi_{2}) \right|_{t=0}, \quad \vec X_{2}(\mat R) = \left.\frac{\d{}}{\d{t}}\mat R(\phi_{1},\Phi+t,\phi_{2}) \right|_{t=0}, \quad \vec X_{3}(\mat R) = \left.\frac{\d{}}{\d{t}}\mat R(\phi_{1},\Phi,\phi_{2}+t) \right|_{t=0}.
\end{equation*}
In this representation, vector fields are typically obtained by differentiating an ODF with respect to the Euler angles $\phi_{1}$, $\Phi$ and $\phi_{2}$.

An explicit description of $\vec X_{1},\vec X_{2},\vec X_{3}$ is given in~\cite{Hielscher2026b}.
From this description, one immediately sees that
$\vec X_{1}(\mat R(\phi_{1},0,\phi_{2})) = \vec X_{3}(\mat R(\phi_{1},0,\phi_{2}))$ and
$\vec X_{1}(\mat R(\phi_{1},\pi,\phi_{2})) = - \vec X_{3}(\mat R(\phi_{1},\pi,\phi_{2}))$
for arbitrary $\phi_{1}$ and $\phi_{2}$.
Hence, this choice of basis degenerates for $\Phi=0$ and $\Phi=\pi$, and for $\Phi\approx0$ or $\Phi\approx\pi$ two of the basis vectors are nearly linearly dependent.

Consequently, the Euler-angle-induced basis, which is often used in the literature, should not be used for stable numerical computations with tangent vector fields.

\paragraph{What is a numerically more favourable basis?}
Two numerically more stable and geometrically more natural choices are obtained directly from the left and right rotation-vector representations introduced above. They are described in detail in~\cite{Hielscher2026b}.

Let $\vec e_1,\vec e_2,\vec e_3$ denote the standard basis of $\IR^{3}$. 
Then, for \(i=1,2,3\), one may choose the three basis vector fields to be given either by
\[ \vec X_i(\mat R) = \mat R\,\widehat{\vec e}_i \qquad\text{or}\qquad  \vec X_i(\mat R) = \widehat{\vec e}_i\,\mat R, \]
which correspond to the right and left tangent vector representations, respectively.
Both choices satisfy properties~1--3 and also the additional orthogonality property~(4).
In particular, they are defined for all orientations $\mat R\in\SO3$, do not degenerate, and are orthonormal at every orientation.

When describing the tangent vectors by the rotation vector, we end up with the vector field representation of description~(c).
More precisely, as the rotation vectors of $\vec X_1,\vec X_2,\vec X_3$ are the standard vectors in $\IR^3$, \cref{eq:LinKombiVectorField} reads as
\[ \vec f_{\star}\colon \SO3 \to \IR^{3}, \quad \vec f_{\star}(\mat R) = \begin{pmatrix} f_{\star,1}(\mat R) \\ f_{\star,2}(\mat R) \\ f_{\star,3}(\mat R) \end{pmatrix}, \qquad \star\in\{L,R\}. \]
We only have to specify whether $\vec f_{\star}(\mat R)$ is understood in the left $(\star=L)$ or in the right $(\star=R)$ representation.
In other words, after choosing either the left or the right representation, a vector field assigns to each orientation the rotation vector of the local orientation change. The functions $f_1,f_2,f_3\colon\SO3\to\IR$ in \cref{eq:LinKombiVectorField} thus serve as the components of the rotation vector.


\section{Symmetries in Vector Fields}\label{sec:symmetries}

In texture analysis, two orientations that differ only by prescribed left and right symmetry operations are usually regarded as physically equivalent, see Morawiec~\cite{Morawiec2004}.
In the classical ODF setting, these are often referred to as specimen and crystal symmetries, depending on the chosen convention.
Here, we use the more general terms left and right symmetries, since in other settings, such as misorientation density functions, both sides may correspond to crystal symmetries.

The symmetries are described by finite subgroups \(S_L,S_R\subset\SO3\) of the orientation space, i.e., the composition of two symmetry operations (rotations) of $S_{L}$ is again an element of $S_{L}$, and analogously for $S_{R}$.
Mathematically, we describe symmetrically equivalent orientations through
\[ \mat R \sim \mat Q_L \, \mat R \, \mat Q_R, \qquad \text{for all } \mat Q_L\in S_L \text{ and } \mat Q_R\in S_R. \]

A scalar orientation-dependent function \(f\colon\SO3\to\IR\), for example an ODF, is called symmetric if it takes the same value at symmetrically equivalent orientations.
For example, quartz has trigonal crystal symmetry with a threefold rotation axis parallel to the crystallographic $c$-axis.
Hence, an ODF of a quartzite takes the same value at any two orientations that differ by a rotation of $120^{\circ}$ about the $c$-axis, i.e., by a right symmetry operation $\mat Q_R$.
Mathematically, we describe this property by
\[ f(\mat R) = f(\mat Q_L \, \mat R \, \mat Q_R), \qquad \text{for all }\mat Q_L\in S_L\text{ and } \mat Q_R\in S_R. \]

\paragraph{How do symmetries act on vector fields?}
For vector fields, the situation is different, since a vector field does not map an orientation to a scalar value, but to a small change of this orientation, i.e., a tangent vector.

To see how such an orientation change behaves under symmetry, we return to the description of small orientation changes from \Cref{sec:TangentVectors}.
Let $\mat R_{0}$ be an orientation and $\mat R_{1}$ a nearby orientation.
Then the relative rotations $\Delta_{L} = \mat R_{1}\mat R_{0}^{\top}$ and $\Delta_{R} = \mat R_{0}^{\top}\mat R_{1}$ describe the change from $\mat R_{0}$ to $\mat R_{1}$, and their rotation vectors $\vec\rho_{L}$ and $\vec\rho_{R}$ approximate the corresponding tangent vector.

We now move this pair as a whole to a symmetrically equivalent position in the orientation space, i.e., we apply the same symmetry operations $\mat Q_{L}$ and $\mat Q_{R}$ to both orientations, which yields $\mat Q_{L}\mat R_{0}\mat Q_{R}$ and $\mat Q_{L}\mat R_{1}\mat Q_{R}$.
This new pair describes the same physical change of orientation as the original one.
Its relative rotations are
\begin{align*}
  \Delta_{L}' &= \mat Q_{L}\mat R_{1}\mat Q_{R}\,\mat Q_{R}^{\top}\mat R_{0}^{\top}\mat Q_{L}^{\top} = \mat Q_{L}\,\Delta_{L}\,\mat Q_{L}^{\top}, \\
  \Delta_{R}' &= \mat Q_{R}^{\top}\mat R_{0}^{\top}\mat Q_{L}^{\top}\,\mat Q_{L}\mat R_{1}\mat Q_{R} = \mat Q_{R}^{\top}\,\Delta_{R}\,\mat Q_{R}.
\end{align*}

The rotation $\Delta_{L}' = \mat Q_{L}\,\Delta_{L}\,\mat Q_{L}^{\top}$ describes the same physical change of orientation as $\Delta_{L}$, only expressed in a coordinate system rotated by $\mat Q_{L}$.
Hence, it is a rotation by the same angle as $\Delta_{L}$, but about a rotation axis that is rotated by $\mat Q_{L}$.
Indeed, if $\vec n$ is the rotation axis of $\Delta_{L}$, i.e., $\Delta_{L}\vec n = \vec n$, then
\[ \Delta_{L}'\,\mat Q_{L}\vec n = \mat Q_{L}\,\Delta_{L}\,\mat Q_{L}^{\top}\mat Q_{L}\,\vec n = \mat Q_{L}\,\Delta_{L}\,\vec n = \mat Q_{L}\,\vec n, \]
and therefore $\mat Q_{L} \vec n$ is the rotation axis of $\Delta_{L}'$.
Analogously, if $\vec n$ denotes the rotation axis of $\Delta_{R}$, the rotation axis of $\Delta_{R}' = \mat Q_{R}^{\top}\,\Delta_{R}\,\mat Q_{R}$ is $\mat Q_{R}^{\top}\vec n$.
Since a rotation vector is the rotation axis scaled by the rotation angle, the left and right rotation vectors transform as
\[ \vec\rho_{L}' = \mat Q_{L}\,\vec\rho_{L} \qquad\text{and}\qquad \vec\rho_{R}' = \mat Q_{R}^{\top}\,\vec\rho_{R}. \]

Accordingly, we call a vector field $\vec f$ symmetric if it maps symmetrically equivalent orientations to correspondingly transformed tangent vectors, i.e.,
\begin{equation}\label{eq:SymmetricVectorField}
  \vec f_{L}(\mat Q_{L}\,\mat R\,\mat Q_{R}) = \mat Q_{L}\,\vec f_{L}(\mat R)
  \qquad\text{or, equivalently,}\qquad
  \vec f_{R}(\mat Q_{L}\,\mat R\,\mat Q_{R}) = \mat Q_{R}^{\top}\,\vec f_{R}(\mat R)
\end{equation}
for all $\mat R\in\SO3$, $\mat Q_{L}\in S_{L}$, and $\mat Q_{R}\in S_{R}$.

Thus, in the left representation, the component functions $f_{L,1},f_{L,2},f_{L,3}$ are invariant under the right symmetry, i.e., $f_{L,i}(\mat R\,\mat Q_{R})=f_{L,i}(\mat R)$, whereas the left symmetry acts on the rotation vectors themselves.
Analogously, in the right representation, the component functions are invariant under the left symmetry, while the right symmetry acts on the rotation vectors.


\section{Harmonic Vector Field Representation}\label{sec:HarmonicVectorFields}

In \Cref{sec:VectorFields}, we have seen that, once the left or the right representation is chosen, an orientation-dependent vector field is described by the rotation vectors of its tangent vectors, i.e., by three scalar component functions
\[ f_{\star,1},f_{\star,2},f_{\star,3} \colon \SO3 \to \IR, \qquad \star\in\{L,R\}. \]
Since these are ordinary orientation-dependent functions, just like an ODF, it is natural to expand them into finite harmonic series of prescribed bandwidth \(N\in\IN\).
More precisely, each component function is represented by its harmonic coefficients and written as
\[ f_{\star,i}(\mat R) = \sum_{n=0}^{N} \sum_{k,k'=-n}^{n} \widehat{(f_{\star,i})}{}_{n}^{k,k'} \, \WignerD{n}{k}{k'}(\mat R), \qquad i=1,2,3.
\]
Here, \(\widehat{(f_{\star,i})}{}_{n}^{k,k'}\) are the harmonic coefficients of the component function \(f_{\star,i}\), and \(\WignerD{n}{k}{k'}\) denote the generalized spherical harmonics, also called Wigner-\(D\) functions.
We refer to a vector field given in this form as a \emph{harmonic vector field}.

\paragraph{Why use the harmonic representation?}
The harmonic representation is one way to describe an arbitrary vector field on the whole orientation space by finitely many coefficients.
Although this is an approximation, its error can be reduced by increasing the bandwidth \(N\).

The more important advantage is that the harmonic representation gives access to fast algorithms.
More precisely, the so-called fast \(\SO3\)-Fourier transform, for which several algorithms are described in~\cite{Risbo1996,Kostelec2008,Potts2009,Hielscher2026}, allows a harmonic vector field to be evaluated very fast at many orientations by fast evaluation of the three scalar component functions.
This can save a considerable amount of computation time.
Consider, for example, the Taylor model in texture evolution.
There, the spin tensor of a crystal with a given orientation is obtained by solving an optimization problem, which is expensive if it has to be done for every orientation in every time step.
Alternatively, one approximates the corresponding vector field once by a harmonic vector field and afterwards evaluates it fast at all required orientations.
In this way, hours of computation time are reduced to seconds.

To this end, it is only necessary to compute the harmonic coefficients of the three component functions.
This can be done efficiently as well.
First, the vector field is sampled on a suitable quadrature grid, i.e., in our example, the spin tensors are computed by the Taylor model only for the orientations of this grid.
Afterwards, a single fast \(\SO3\)-Fourier transform yields the harmonic coefficients.

\paragraph{Left or right representation?}
Before computing the harmonic representation, one has to choose whether the vector field is represented in left or right rotation-vector coordinates.

In general, the corresponding component functions and therefore also their harmonic coefficients are different in the two representations.
However, by~\cite[Lemma~5]{Hielscher2026b}, the change between the two representations almost preserves the harmonic bandwidth.
Moreover, as shown in~\cite{Hielscher2026b}, this change can be computed directly in the frequency domain.
This means that the harmonic coefficients \(\widehat{(f_{R,i})}{}_{n}^{k,k'}\) can be computed easily from the coefficients \(\widehat{(f_{L,i})}{}_{n}^{k,k'}\), and vice versa.
In the MTEX implementation, however, this transformation is computed by discrete fast \(\SO3\)-Fourier algorithms, since practical computations showed this to be significantly faster.

Since both representations are equally suitable for harmonic approximation and can be converted into each other, the choice of representation should be guided by symmetry.
As discussed in \Cref{sec:symmetries}, the component functions in the left representation are invariant under the right symmetry group, whereas the component functions in the right representation are invariant under the left symmetry group.
Since symmetry of orientation-dependent functions translates into symmetry relations for their harmonic coefficients, as described in~\cite{Bunge1982,Hielscher2026}, which can be exploited to speed up the fast algorithms, we should choose the representation in which the larger symmetry group acts directly on the component functions.
Consequently, if $\abs{S_L} > \abs{S_R}$, we use the right representation, and otherwise the left representation.

\section{Computation of Gradient, Divergence, and Curl}\label{sec:GradientDivCurl}

A further advantage of the harmonic representation of vector fields is that it extends the harmonic framework of orientation-dependent functions.
In particular, differential operators such as the gradient, the divergence, and the curl can be computed directly from the harmonic coefficients, see~\cite{Hielscher2026b}, so that no numerical differentiation, e.g., by finite differences, is required.
In the following, we summarize the corresponding formulas.

\paragraph{The gradient:}
The gradient of an $N$-band-limited scalar orientation-dependent function
\[ f(\mat R) = \sum_{n=0}^{N} \sum_{k,k'=-n}^{n} \fhat{n}{k}{k'} \, \WignerD{n}{k}{k'}(\mat R) \]
at an orientation $\mat R$ is the tangent vector at $\mat R$ that points in the direction of steepest increase of $f$ and whose length equals the rate of this increase.
Since this tangent vector varies smoothly with $\mat R$, the gradient $\nabla f$ is a vector field, which we describe, as before, by its rotation vectors $(\nabla f)_{L}$ or $(\nabla f)_{R}$.

By~\cite[Corollary~9]{Hielscher2026b}, the component functions of $(\nabla f)_{L}$ and $(\nabla f)_{R}$ are again $N$-band-limited.
Hence, the gradient of $f$ is a harmonic vector field of the same bandwidth $N$ as $f$ itself.
Moreover, by~\cite[Corollary~8]{Hielscher2026b}, its harmonic coefficients can be computed directly from the harmonic coefficients $\fhat{n}{k}{k'}$ of $f$.
In the left representation, we obtain the harmonic vector field
\begin{equation*}
  (\nabla f)_{L} (\mat R) =
  \begin{pmatrix}
    \ds\sum_{n=0}^{N}\sum_{k,k'=-n}^n \left(-\i \, c_n^{-k} \, \fhat{n}{k-1}{k'} - \i \, c_n^k \, \fhat{n}{k+1}{k'}\right) \, \WignerD{n}{k}{k'}(\mat R) \\
    \ds\sum_{n=0}^{N}\sum_{k,k'=-n}^n \left(c_n^{-k} \, \fhat{n}{k-1}{k'} - c_n^k \, \fhat{n}{k+1}{k'}\right) \, \WignerD{n}{k}{k'}(\mat R) \\
    \ds\sum_{n=0}^{N}\sum_{k,k'=-n}^n \left(-\i\,k\,\fhat{n}{k}{k'}\right) \, \WignerD{n}{k}{k'}(\mat R)
  \end{pmatrix},
\end{equation*}
with $c_n^m = \frac12\sqrt{(n-m)(n+m+1)}$.
The terms in brackets are the harmonic coefficients of the three component functions, i.e., each of them is a weighted sum of at most two coefficients of $f$ with the same or neighboring indices.
In the right representation, we obtain
\begin{equation*}
  (\nabla f)_{R} (\mat R) =
  \begin{pmatrix}
    \ds\sum_{n=0}^{N}\sum_{k,k'=-n}^n \left(-\i \, c_n^{-k'} \, \fhat{n}{k}{k'-1} - \i \, c_n^{k'} \, \fhat{n}{k}{k'+1}\right) \, \WignerD{n}{k}{k'}(\mat R) \\
    \ds\sum_{n=0}^{N}\sum_{k,k'=-n}^n \left(-c_n^{-k'} \, \fhat{n}{k}{k'-1} + c_n^{k'} \, \fhat{n}{k}{k'+1}\right) \, \WignerD{n}{k}{k'}(\mat R) \\
    \ds\sum_{n=0}^{N}\sum_{k,k'=-n}^n \left(-\i\,k'\,\fhat{n}{k}{k'}\right) \, \WignerD{n}{k}{k'}(\mat R)
  \end{pmatrix}.
\end{equation*}
We use the convention that harmonic coefficients with indices outside the admissible range are zero, i.e.,
\[ \fhat{n}{k}{k'}=0 \qquad \text{whenever } |k|>n \text{ or } |k'|>n. \]

\paragraph{The divergence:}
The divergence of a vector field describes the local flux balance of the field.
It is a scalar orientation-dependent function: positive values correspond to local sources, whereas negative values correspond to local sinks.
The divergence can be computed directly in the frequency domain by applying the gradient formulas above to the scalar component functions.
More precisely, by~\cite[Lemma~10]{Hielscher2026b},
\begin{align*}
  \div\vec f &= \braces[\big]{\grad f_{L,1}}_{L,1} + \braces[\big]{\grad f_{L,2}}_{L,2} + \braces[\big]{\grad f_{L,3}}_{L,3} \\
   &= \braces[\big]{\grad f_{R,1}}_{R,1} + \braces[\big]{\grad f_{R,2}}_{R,2} + \braces[\big]{\grad f_{R,3}}_{R,3}.
\end{align*}
Here, \(\braces[\big]{\grad f_{L,1}}_{L,1}\) denotes the first component of the gradient of the scalar function \(f_{L,1}\), computed in the left rotation-vector representation.

An important application of the divergence is texture evolution under plastic deformation.
There, the evolution of the ODF $f$ is governed by the transport equation, cf.~\cite{Clement1979,Morawiec2004},
\[ \partial_t f(t,\mat R) = - \div\big(f(t,\mat R) \, \vec\Omega(t,\mat R)\big), \]
where $\vec\Omega$ denotes the spin field, i.e., the vector field that assigns to every crystal orientation the rotation vector of its lattice rotation.

\paragraph{The curl:}
The curl of a vector field is again a vector field.
While the divergence describes local sources and sinks, the curl describes the local circulation of the field, i.e., how strongly the field swirls around a given orientation.

A vector field has vanishing curl if and only if it is the gradient of a scalar function $g\colon\SO3\to\IR$.
In particular, $\curl(\nabla g) = 0$ for every such function~$g$.
Hence, the curl allows one to check whether a given vector field is a gradient field or not.

By~\cite[Lemma~10]{Hielscher2026b}, the curl is represented in left rotation-vector coordinates by
\begin{align*}
  (\curl\vec f)_{L} =
  \begin{pmatrix*}
    \braces[\big]{\grad f_{L,2}}_{L,3} - \braces[\big]{\grad f_{L,3}}_{L,2} - f_{L,1} \\
    \braces[\big]{\grad f_{L,3}}_{L,1} - \braces[\big]{\grad f_{L,1}}_{L,3} - f_{L,2} \\
    \braces[\big]{\grad f_{L,1}}_{L,2} - \braces[\big]{\grad f_{L,2}}_{L,1} - f_{L,3}
  \end{pmatrix*}.
\end{align*}
In right rotation-vector coordinates, it is represented by
\begin{align*}
  (\curl\vec f)_{R} =
  \begin{pmatrix*}
    \braces[\big]{\grad f_{R,2}}_{R,3} - \braces[\big]{\grad f_{R,3}}_{R,2} + f_{R,1} \\
    \braces[\big]{\grad f_{R,3}}_{R,1} - \braces[\big]{\grad f_{R,1}}_{R,3} + f_{R,2} \\
    \braces[\big]{\grad f_{R,1}}_{R,2} - \braces[\big]{\grad f_{R,2}}_{R,1} + f_{R,3}
  \end{pmatrix*}.
\end{align*}


\subsection*{Use of AI}
AI-based tools were used to check spelling and grammar and to improve the clarity of the English. They were not used to generate any scientific content.

\printbibliography

\end{document}